# Combining Homotopy Methods and Numerical Optimal Control to Solve Motion Planning Problems

Kristoffer Bergman and Daniel Axehill

*Abstract*— This paper presents a systematic approach for computing local solutions to motion planning problems in non-convex environments using numerical optimal control techniques. It extends the range of use of state-of-the-art numerical optimal control tools to problem classes where these tools have previously not been applicable. Today these problems are typically solved using motion planners based on randomized or graph search. The general principle is to define a homotopy that perturbs, or preferably relaxes, the original problem to an easily solved problem. By combining a Sequential Quadratic Programming (SQP) method with a homotopy approach that gradually transforms the problem from a relaxed one to the original one, practically relevant locally optimal solutions to the motion planning problem can be computed. The approach is demonstrated in motion planning problems in challenging 2D and 3D environments, where the presented method significantly outperforms a state-of-the-art open-source optimizing sampled-based planner commonly used as benchmark.

## I. INTRODUCTION

### A. Background

Depending on the field of subject and the applications, the meaning of the term "optimal motion planning" may vary [1]. In this paper, the term will be used to denote the problem of finding a feasible and locally optimal trajectory for a dynamical system from an initial state to a goal state in an environment that includes obstacles. Hence, this is a quite broad notion which includes several applications such as finding trajectories for autonomous cars, unmanned aerial vehicles, autonomous underwater vehicles and trajectories for robot manipulators.

In the literature, there are different methods that can be used to solve motion planning problems and produce feasible or optimal trajectories. For more information and details on some of these methods, the reader is referred to standard literature such as [1].

In general, finding a feasible trajectory for a system with nonlinear dynamics in complex non-convex environments, such as bugtraps and mazes, is a challenging problem. Most of the work in this area has been focused on sampled-based and graph search based methods, which are able to handle sampled maps and challenging nonlinear dynamics [1]. However, these methods suffer from the *curse of dimensionality* as the dimension of the search space to explore grows. Furthermore, lattice planners, a commonly used graph search method, typically rely on a preprocessing step offline where so-called motion primitives are precomputed [1]. This step limits the flexibility of the method when applied to problems with, *e.g.*, dynamics that change over time.

Another approach is to use numerical optimal control methods to solve the motion planning problem. However, the difficulty with these methods is that they are sensitive to initialization, especially in non-convex environments [2]. As a result, these methods are considered less suitable and are rarely used for these types of problems.

### B. Related Work

In recent years, the interest within research on sampling-based motion planning methods with focus on optimality, such as Rapidly-exploring Random Trees (RRT*), Probabilistic Roadmaps (PRM*) [3] and Stable Sparse RRT (SST and SST*) [4], has increased. The original versions of RRT* and PRM* can not handle nonholonomic dynamical constraints. An extension of the RRT* algorithm to such problem classes is presented in [5]. However, this algorithm is dependent on a steering function, which in the case of exact steering for a general nonholonomic dynamic system corresponds to solving a boundary value problem to connect samples in the tree. Another approach is to use RRT to sample the input to a closed-loop system with an underlying controller [6].

One recent method of graph search based planning combined with numerical optimal control is [7], where a state lattice is generated for a general 2-trailer system using numerical optimization, which can be combined with any graph-based search algorithm, such as A*. The method is computationally efficient since only a graph search needs to be performed online. However, the solution is only resolution complete, and it requires a computationally demanding pre-processing step where motion primitives are computed and the search tree is prepared.

Using optimal control to solve motion planning problems has been investigated by many researchers. In the early paper of [8] the authors propose an optimal control approach to compute a minimum time solution to a motion planning problem. However, the desired geometric path is assumed given. In [9], a motion planning problem for a robot in an environment with obstacles is considered. Here, both the obstacle constraints and the system dynamical equations are treated simultaneously by incorporating methods from numerical optimal control to the planning problem.

One recent publication where the use of newly developed tools for numerical optimal control are used to study vehicle-maneuver and trajectory generation is [10]. In that publication three different maneuvers are considered for models of varying complexity, and the optimal control solutions are

This work was partially supported by the Wallenberg Autonomous Systems and Software Program (WASP).

K. Bergman and D. Axehill are with the Division of Automatic Control, Linköping University, Sweden {kristoffer.bergman, daniel.axehill}@liu.se

presented and analyzed. The result indicates that it is often non-trivial to find an initial guess that solves the problem.

Homotopy, or continuation, based methods have been used to solve manipulation planning problems [11], where the constraints are relaxed and iteratively introduced to obtain the original problem formulation. The relaxed subproblems are solved using a geometric planner (*i.e.* nonholonomic dynamical systems are not considered) based on dynamic programming techniques. Another example is [12], where planning for a nonholonomic mobile robot is carried out by first using a geometric planner to satisfy all geometric constraints, and then introduce the nonholonomic constraints iteratively.

Homotopy methods have also been used to solve two-point boundary value problems [13], [14], where the homotopy method has been combined with iterative search methods to compute solutions to problems with nonlinear dynamics when the initial guess is far from the optimal solution. However, the case with nonlinear inequality constraints representing non-convex obstacles has not previously been considered.

*C. Main Contribution*

In this work, it is shown how a state-of-the-art method from numerical optimal control can be combined with a homotopy method into a method that solves motion planning problems with challenging geometry of the feasible set. The possibility to combine the homotopy method and the SQP method has previously been considered for generic optimization in, *e.g.*, [15]. However, the combination has not been used for motion planning problems before. Furthermore, the method is combined with a novel structure-exploiting constraint parameterization introduced in this paper. As a result this is, to the best of the authors' knowledge, the first method that in a structured way enables the use of the powerful methods known from numerical optimal control on challenging motion planning problems. It is illustrated that the proposed method completely outperforms an optimizing sampled-based planner in terms of computational effort and trajectory quality on examples where a traditional state-of-the-art numerical optimal control method fails.

## II. OPTIMAL CONTROL PRELIMINARIES

The traditional optimal control problem is to compute an optimal control input to a dynamical system such that an objective function is minimized, subject to initial and terminal constraints [16]. By also enforcing the solution to satisfy certain inequality constraints, the problem will be denoted a *constrained continuous-time optimal control* (CCTOC) problem. It is defined in its general form as

$$\begin{aligned}
\underset{x(t),u(t),t_f}{\text{minimize}} \quad & \Phi(x(t_f)) + \int_{t_0}^{t_f} f_0(x(t), u(t)) dt \\
\text{subject to} \quad & x(t_0) = x_{t_0}, \\
& \dot{x}(t) = f(x(t), u(t)), \, \forall t \in [t_0, t_f] \\
& g(x(t), u(t)) \preceq 0, \, \forall t \in [t_0, t_f] \\
& \Psi(x(t_f)) = x_{t_f},
\end{aligned} \quad (1)$$

where $x_{t_0}$ and $x_{t_f}$ are the initial and final states, $x(t) \in \mathbb{R}^n$ are the states, $u(t) \in \mathbb{R}^m$ are the control inputs, $f : \mathbb{R}^n \times \mathbb{R}^m \to \mathbb{R}^n$ describes the system dynamics, $g : \mathbb{R}^n \times \mathbb{R}^m \to \mathbb{R}^p$ are the inequality constraint functions that determine the feasible region, $\Phi : \mathbb{R}^n \to \mathbb{R}$ and $f_0 : \mathbb{R}^n \times \mathbb{R}^m \to \mathbb{R}$ form the objective function and $\Psi : \mathbb{R}^n \to \mathbb{R}^n$ is the terminal constraint function. The final time $t_f$ may, or may not, be an optimization variable [16].

*A. Solving the Optimal Control Problem*

The CCTOC problem (1) can be solved using several different approaches, see *e.g.* [17], [18], [2]. In this paper, we are using a direct method, in which the constrained optimal control problem is reformulated into a nonlinear programming (NLP) problem. The main reasons for this choice are that indirect formulations are often hard to solve in practice [2], and that there exist several good software packages for solving the reformulated NLP problem with these methods, such as ACADO [19], IPOPT [20], SNOPT [21] and several more.

In general, the solution to an NLP problem is found using iterative schemes [22]. The very basic idea of an iterative solver is presented in Algorithm 1. The major difference between different iterative solvers is how to compute the update direction and step size. Different solvers also use different means to avoid getting stuck in local optimal solutions: some do not care at all about this, while there exist methods like, *e.g.*, genetic algorithms (GAs) where an attempt to avoid these are made [23]. An important aspect to get good solutions (*i.e.* local solutions of practical relevance) is the choice of the initial iterate $x_0$ in Algorithm 1, and for some methods there exist structured ways or good heuristics to do this. For more information on numerical solvers for NLP problems, the reader is referred to standard literature such as [22].

---

**Algorithm 1** Outline of iterative NLP solver

---

1: Given an initial iterate $x_0$
2: $k \leftarrow 0$
3: **while** $x_k$ not locally optimal **do**
4:     Compute an update direction $p_k$ and step size $\alpha_k$
5:     $x_{k+1} \leftarrow x_k + \alpha_k p_k$
6: **end while**

---

In this work, we have used ACADO for solving the CCTOC problem. ACADO uses an SQP method combined with support for different types of discretization methods and integrators to solve the CCTOC problem (see [19] for further details). In order to obtain a finite dimensional NLP problem, we have used the direct multiple shooting discretization method, combined with the classical Runge-Kutta integrator [2].

*B. Nonlinear Programming Formulation*

When the discretization method for the CCTOC problem is decided, it can be formulated as an NLP problem. Since we are using a multiple shooting approach, the problem is formulated as a collection of $N$ phases. If we assume piecewise constant control signals $u_i$, each interval can be solved numerically, starting from an artificial initial state value $s_i$

$$\dot{x}_i(t) = f(x_i(t), u_i), \quad t \in [t_i, t_{i+1}], \\ x_i(t_i) = s_i. \quad (2)$$

This produces trajectory sequences $x_i(t; s_i, u_i)$ for each phase that are determined by the initial state values $s_i$ and the constant control signal $u_i$. The NLP decision variables are then $[s_0, u_0, s_1, u_1, \ldots, s_N]$ [24], which gives the following NLP formulation of the CCTOC problem, assuming that the cost function for the final state is zero and that the terminal constraint function is the identity function

$$\begin{aligned} \underset{\mathbf{s}, \mathbf{u}}{\text{minimize}} \quad & \sum_{i=0}^{N-1} l_i(s_i, u_i) = L(\mathbf{s}, \mathbf{u}) \\ \text{subject to} \quad & s_0 = x_{t_0}, \\ & s_{i+1} = x_i(t_{i+1}; s_i, u_i), \ \forall i \in \mathbb{Z}_{0, N-1} \\ & g(s_i, u_i) \preceq 0, \ \forall i \in \mathbb{Z}_{0, N} \\ & s_N = x_{t_f}. \end{aligned} \quad (3)$$

Here, $\mathbf{s}, \mathbf{u}$ are the vectors containing all the discretized states and control signals, and the constraints $s_{i+1} = x_i(t_{i+1}; s_i, u_i)$ are the so-called *continuity* constraints that connect the phases. The objective function is a sum of the numerically computed integrals for each phase [24].

*C. Sequential Quadratic Programming*

Once the CCTOC problem has been reformulated as an NLP problem, it is possible to solve it using, *e.g.*, SQP. If we define the notation $\mathbf{w} = [s_0, u_0, s_1, u_1, \ldots, s_N]$ for the optimization variables in (3), an SQP method computes a search direction $\mathbf{p}$ by approximating the NLP by a quadratic program (QP) at iterate $(\mathbf{w}_k, \lambda_k)$, where $\lambda_k$ are the Lagrange multipliers from the solution in the last iterate, by solving

$$\begin{aligned} \underset{\mathbf{p}}{\text{minimize}} \quad & \frac{1}{2} \mathbf{p}^T \mathbf{B}_k \mathbf{p} + \nabla L(\mathbf{w}_k)^T \mathbf{p} \\ \text{subject to} \quad & \nabla F(\mathbf{w}_k)^T \mathbf{p} + F(\mathbf{w}_k) = 0, \\ & \nabla G(\mathbf{w}_k)^T \mathbf{p} + G(\mathbf{w}_k) \preceq 0. \end{aligned} \quad (4)$$

Here, $F$ contains all the discrete time dynamic system constraints as well as initial and final constraints and $G$ all the inequality constraint functions. $\mathbf{B}_k$ denotes the approximation of the Hessian of the Lagrangian to the NLP, $\nabla F$, $\nabla G$ the Jacobian matrices of $F$ and $G$ respectively, and $\nabla L$ is the gradient of the objective function. For more details, see [22].

Due to the equivalence between SQP and Newton's method, the search direction $\mathbf{p}$ can either be defined as the solution to the quadratic program (4), or as the search direction generated by Newton's method for the corresponding NLP problem (3) for any fixed working set (*i.e.* the set of inequality constraints that holds with equality). In [22] it is stated that the SQP method eventually will find the optimal active set during the SQP iterations under certain mild conditions. Once that is done, what remains is basically Newton's method on a nonlinear system of equations.

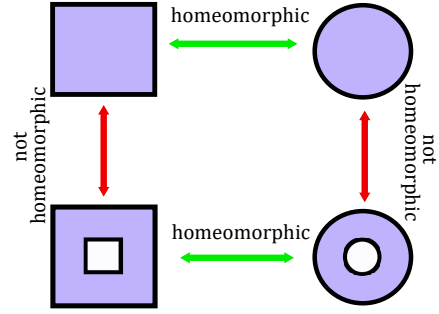

Fig. 1: An example illustrating homeomorphic and non-homeomorphic transformations of objects in 2D. The transformations from the top objects to the bottom objects are not homeomorphic due to the hole in the middle.

### III. CONTINUATION/HOMOTOPY METHODS

Since solvers for a system of nonlinear equations are sensitive to the choice of initial iterate $x_0$, they are in danger of converging to a local minimum which does not coincide with a global minimum, or not converging at all, for difficult problems. Convergence to a globally optimal solution might be too much to ask for, but what is desired in practice is convergence to a local solution of practical relevance. One way to deal with this is to use so-called *homotopy*, or *continuation* methods [25]. Instead of solving the original problem directly, one can set up an easy problem for which the solution is more or less obvious, and then successively transform the problem towards the original and more difficult one.

*Example 1:* If the problem is to find a root for the nonlinear equation $F(x) = 0$ ($F : \mathbb{R}^n \to \mathbb{R}$), the homotopy method instead focuses on finding a root to the *homotopy map*

$$H(x, \gamma) = \gamma F(x) + (1 - \gamma) G(x) \quad (5)$$

where $\gamma \in [0, 1]$, from now on denoted *homotopy parameter*, and $G : \mathbb{R}^n \to \mathbb{R}$. The basic idea is to solve a sequence of $n+1$ problems $H(x, \gamma) = 0$, for $\gamma = \gamma_0 < \gamma_1 < \ldots < \gamma_n = 1$. If the roots are relatively easy to find for $G(x)$, it should be easy to find a first iterate $x_0$. Then for any problem $k \in [0, n]$ along the homotopy path, given a starting point $x_k$ that is the approximate solution from $H(x, \gamma_k) = 0$, the next approximate solution $x_{k+1}$ is calculated by solving $H(x, \gamma_{k+1}) = 0$. This can be done using an iterative method such as Newton's method, with $x_k$ as initial iterate. If the difference between $\gamma_k$ and $\gamma_{k+1}$ is small enough, the method will converge [25].

If we separate the inequality constraints $g$ into two parts

- $g_C : \mathbb{R}^n \times \mathbb{R}^m \to \mathbb{R}^r$, that contains all the convex inequality constraint functions and
- $g_H : \mathbb{R}^n \times \mathbb{R}^m \to \mathbb{R}^q$, $r + q = p$, that contains all the non-convex inequality constraint functions,

it is possible to replace the non-convex inequality constraints (in our case representing the obstacles) with the homotopy map $h : \mathbb{R}^n \times \mathbb{R}^m \times \mathbb{R} \to \mathbb{R}^q$, which represents a continuous transformation of $g_H$ with respect to the homotopy parameter $\gamma$. We require that $h$ is a *homeomorphism*[1] with respect to the inequality constraint functions.

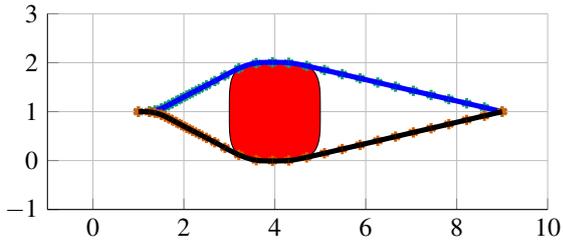

Fig. 2: An example illustrating an environment with two homotopy classes. The trajectories can not be continuously transformed into each other due to the obstacle.

*Definition 1:* A function $f : X \to Y$ between two topological spaces $X$ and $Y$ is called a homeomorphism if
1) $f$ is a one-to-one and onto mapping (bijection), and
2) both $f$ and $f^{-1}$ are continuous [1].

An illustrating example of homeomorphic transformations can be found in Figure 1. The reason to require that $h$ is a homeomorphism is easily motivated in the context of an CCTOC problem.

*Definition 2:* Two trajectories are said to be in the same *homotopy class* if they can be continuously transformed from one to another without intersecting any obstacles, while keeping their endpoints fixed [28].

An example of an environment with two homotopy classes is shown in Figure 2: one class above the obstacle and one class below. Hence, the homeomorphism requirement on $h$ prevents the introduction of new, and removal of existing, homotopy classes, while continuously perturbing the available free space in each of the existing homotopy classes when $\gamma$ is updated.

Now, using the homotopy map, we can write the NLP in (3) as

$$\begin{aligned}
\underset{\mathbf{s},\mathbf{u}}{\text{minimize}} \quad & L(\mathbf{s},\mathbf{u}) \\
\text{subject to} \quad & s_0 = x_{t_0} \\
& s_{i+1} = x_i(t_{i+1}; s_i, u_i), \forall i \in \mathbb{Z}_{0,N-1} \\
& g_C(s_i, u_i) \preceq 0, \forall i \in \mathbb{Z}_{0,N} \\
& h(s_i, u_i, \gamma) \preceq 0, \forall i \in \mathbb{Z}_{0,N} \\
& s_N = x_{t_f}.
\end{aligned} \quad (6)$$

With this formulation, we can solve the problem iteratively by first solving a problem with simple or even convex inequality constraint functions, and then gradually perturb the non-convex inequality constraints by increasing the $\gamma$-parameter. This method shares similarities with the *barrier method*, where the inequality constraints are relaxed and added to the objective function [26]. Then, a sequence of problems is solved, starting from a simpler problem and then approaching the solution to the original one.

---

[1] The method can be generalized to handle a sequence of homeomorphic mappings.

One important special type of perturbation structure is when the optimization problem is gradually relaxed (in optimization sense) as the homotopy parameter is decreased from 1 to 0. In this case, the objective function value for iterate $[\mathbf{s}_k, \mathbf{u}_k, \gamma_k]$ will always be greater than or equal to the objective value for the previous iterate $[\mathbf{s}_{k-1}, \mathbf{u}_{k-1}, \gamma_{k-1}]$ as the problem is solved for a strictly increasing $\gamma$, *i.e.*

$$L^*(\mathbf{s}_k, \mathbf{u}_k) \geq L^*(\mathbf{s}_{k-1}, \mathbf{u}_{k-1}), \quad \gamma_k > \gamma_{k-1}, \quad (7)$$

as the homotopy parameter is increased from 0 to 1.

## IV. SEQUENTIAL HOMOTOPY QUADRATIC PROGRAMMING

In this section, a prototype algorithm forming the backbone of the main contribution in this paper will be explained. The SQP subproblem (4) can be extended to also include the homotopy map

$$\begin{aligned}
\underset{\mathbf{p}}{\text{minimize}} \quad & \frac{1}{2}\mathbf{p}^T \mathbf{B}_k \mathbf{p} + \nabla L(\mathbf{w}_k)^T \mathbf{p} \\
\text{subject to} \quad & \nabla F(\mathbf{w}_k)^T \mathbf{p} + F(\mathbf{w}_k) = 0, \\
& \nabla G_c(\mathbf{w}_k)^T \mathbf{p} + G_c(\mathbf{w}_k) \preceq 0, \\
& \nabla H(\mathbf{w}_k, \gamma)^T \mathbf{p} + H(\mathbf{w}_k, \gamma) \preceq 0,
\end{aligned} \quad (8)$$

where the inequalities have been divided into $G_c$, which contains all the convex inequality constraint functions, and $H$ which contains all the homotopy mappings. We observe that for a fixed $\gamma$, the problem can be solved, as usual, by a standard SQP method to local optimality. Furthermore, for a fixed working set, the Karush-Kuhn-Tucker (KKT) system is a system of nonlinear equations for which a standard homotopy approach as the one described in [25] can be used. Similarly to [15], we now propose to combine these two methods into one. We call the resulting method Sequential Homotopy Quadratic Programming (SHQP), which is described in Algorithm 2, using notations from (8). It is inspired by Algorithm 18.3 in [22], which describes an SQP method using line search combined with a merit function to globalize the search for a solution. For more details regarding SQP, the reader is referred to [22].

The outer iteration in the algorithm corresponds to the homotopy mapping, while the inner iteration corresponds to the standard SQP algorithm. For every value of $\gamma$, the SQP steps are iterated until the termination test is completed, or until a maximum number of $N_{\text{SQP}}$ SQP iterations has been performed. The termination test can vary with $\gamma$. For example, it might not be necessary to solve the intermediate problems to convergence, but instead accept a lower solution accuracy in the termination test. In the last step however, we would like the termination test to guarantee convergence. The same approach is used for the *inexact centering* in the barrier method [26], where the barrier parameter (homotopy parameter) is updated before Newton's method has fully converged.

An important question to consider is how to update $\gamma$ in every step, *i.e.*, how to choose the step size $\Delta_\gamma$, to obtain the fastest possible perturbation. We want to use as large steps as possible, such that the problem is solved with few homotopy iterations. However, if $\Delta_\gamma$ is too large, the initialization from the last step could be far from the solution to the updated problem. This will require many SQP iterations before the termination test is satisfied in Algorithm 2, or it might not converge at all. An illustrating example of this can be seen in Figure 3. This trade-off is shared both with standard homotopy methods as well as barrier methods [25], [26].

**Algorithm 2** Sequential Homotopy Quadratic Programming
1: Choose line search and merit function parameters [22]
2: $\gamma \leftarrow \gamma_{\text{init}}$, $\mathbf{w}_0 \leftarrow \mathbf{w}_{\text{init}}$, $\lambda_0 \leftarrow \lambda_{\text{init}}$
3: Choose an initial positive definite Hessian approx. $\mathbf{B}_0$
4: Evaluate $L(\mathbf{w}_0), \nabla L(\mathbf{w}_0), F(\mathbf{w}_0), \nabla F(\mathbf{w}_0), G_c(\mathbf{w}_0),$
   $\nabla G_c(\mathbf{w}_0)$
5: **while** $\gamma \leq 1$ **do**
6:    Evaluate $H(\mathbf{w}_0, \gamma), \nabla H(\mathbf{w}_0, \gamma)$
7:    **for** $k = 0, 1 \ldots N_{\text{SQP}}$ **do**
8:      **if** Termination test satisfied **then**
9:        STOP
10:      **end if**
11:      Compute $\mathbf{p}_k$ by solving (8) at $\mathbf{w}_k, \lambda_k, \gamma$
12:      Calculate step size parameter $\alpha_k \in [\alpha_{\min}, 1]$ [22]
13:      $\mathbf{w}_{k+1} \leftarrow \mathbf{w}_k + \alpha_k \mathbf{p}_k$
14:      Evaluate $L(\mathbf{w}_{k+1}), \nabla L(\mathbf{w}_{k+1}), F(\mathbf{w}_{k+1}), \nabla F(\mathbf{w}_{k+1}),$
   $G_c(\mathbf{w}_{k+1}), \nabla G_c(\mathbf{w}_{k+1}), H(\mathbf{w}_{k+1}, \gamma), \nabla H(\mathbf{w}_{k+1}, \gamma)$
15:      Compute new Lagrange multipliers $\lambda_{k+1}$ [22]
16:      Update Hessian using a quasi-Newton formula [22]
17:    **end for**
18:    $\gamma \leftarrow \gamma + \Delta_\gamma$
19:    $\mathbf{w}_0 \leftarrow \mathbf{w}_{k+1}, \lambda_0 \leftarrow \lambda_{k+1}, \mathbf{B}_0 \leftarrow \mathbf{B}_{k+1}$
20:    $L(\mathbf{w}_0), F(\mathbf{w}_0), G_c(\mathbf{w}_0) \leftarrow L(\mathbf{w}_{k+1}), F(\mathbf{w}_{k+1}), G_c(\mathbf{w}_{k+1})$
21:    $\nabla L(\mathbf{w}_0), \nabla F(\mathbf{w}_0), \nabla G_c(\mathbf{w}_0) \leftarrow \nabla L(\mathbf{w}_{k+1}), \nabla F(\mathbf{w}_{k+1}),$
   $\nabla G_c(\mathbf{w}_{k+1})$
22: **end while**

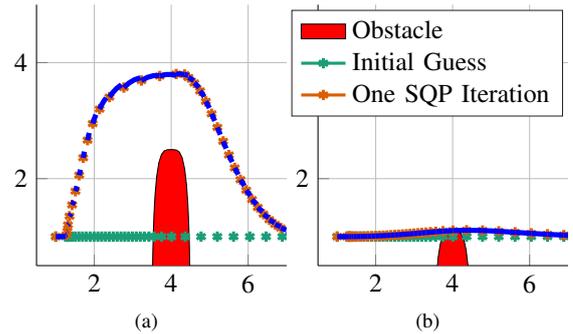

Fig. 3: An example of a large step (a), and a small step (b) in $\gamma$, when an obstacle is perturbed by pushing it into the free space. The blue line represents the simulated trajectory from every multiple shooting point. The small step has already converged after one SQP iteration, while the large step is far from convergence.

The approach can be seen as a predictor-corrector method, with a trivial predictor step that only updates the homotopy parameter. Then, the SQP iterations are considered as the corrector step. If first order information regarding the solutions is available, the predictor step can simultaneously update both the homotopy parameter and the optimization variables, similarly to what is described in [15].

This method can be used to generate trajectories offline or, because of the good warm-start properties of the SQP method in general, be used in real-time online using constantly updated information about the state and the environment. Hence, the proposed method has the potential to be used as a reactive planner.

## V. OBSTACLE CLASSIFICATION

The total number of SQP iterations in Algorithm 2 depends on several factors, such as $\Delta_\gamma$ and the nonlinear dynamics of the system. Other important factors are how large and what type of obstacles that are included in the CCTOC problem. We have categorized the obstacles in two main classes based on topological properties within the obstacle region. We denote these classes *Connected* and *Disconnected*.

The *Connected* obstacles are, in the original problem, connected with the complement to the feasible set of one convex inequality constraint. This specifies in which direction the *Connected* obstacle will grow during the homotopy phase, such that the space representing the feasible area does not introduce new homotopy classes. The name refers to the fact that in this case the infeasible area is topologically connected.

The *Disconnected* obstacles are obstacles *not* connected with the convex obstacle region. In this case, the homotopy map can, *e.g.*, be introduced with a point representing the center of the obstacle. Introducing a Disconnected obstacle will increase the number of homotopy classes for the trajectory planning problem.

Both these types of obstacles can be further divided into *Single* and *Chained* variants of the main classes. A *Single*-type obstacle is represented only by one constraint, while the *Chain*-type obstacles are represented by more than one constraint which are geometrically connected in a chain. The classification of obstacles will be of importance when choosing a homotopy map that satisfies the homeomorphism requirement. Example of these types of obstacles can be found in Figure 4. The chosen homotopy maps for (a)-(c) gradually introduces the obstacles from the outside of the free space. In (b), the long Chain-type obstacle is gradually extended like a "snake". In (d), the radius of the torus is large in the beginning and gradually reduced.

## VI. EXAMPLES OF OBSTACLE REPRESENTATIONS

Algorithm 2 accepts all sorts of obstacle constraint representations that can be handled by SQP solvers, which means that the constraint function should be smooth. One example of practical relevance is to represent them as super-ellipsoids

$$\left(\frac{x - x_c}{r_x}\right)^k + \left(\frac{y - y_c}{r_y}\right)^k + \left(\frac{z - z_c}{r_z}\right)^k - 1 \geq 0, \quad (9)$$

for $k = 2N$, $N \in \mathbb{Z}^+$. $[x_c, y_c, z_c]$ represents the center-point of the obstacle, while $[r_x, r_y, r_z]$ defines the three radii for the super-ellipsoid in each axis. Higher values of $k$ produces obstacles with sharper edges. However, an increasing $k$ will cause the linear approximation of the constraints used in the SQP iterations to be more numerically sensitive. This will in turn lead to a slower convergence. Some examples of super-ellipsoid obstacles with $k = 4$ can be found in Figure 4.

*Example 2:* Assume that the obstacle in Figure 3 is defined by the inequality

$$(x-4)^4 + \left(\frac{y}{6}\right)^4 - 1 \geq 0. \quad (10)$$

This inequality can be replaced by the homotopy map

$$(x-4)^4 + \left(\frac{y-6(1-\gamma)}{6}\right)^4 - 1 \geq 0, \quad (11)$$

which will continuously move the obstacle into the free space when $\gamma$ is increased from 0 to 1. For $\gamma = 1$, we have the original formulation (10).

A second example we have used is to model an obstacle as a torus (illustrated in Figure 4 (d))

$$\left(R - \sqrt{(x-x_c)^2 + (y-y_c)^2}\right)^2 + (z-z_c)^2 - r^2 \geq 0). \quad (12)$$

Here, $R$ defines the distance from the center of the torus to the center of the tube, and $r$ is the radius of the tube.

## VII. NUMERICAL RESULTS

In this section the numerical results from four different trajectory planning scenarios are presented. To produce the results, we have emulated Algorithm 2 for the sake of proof-of-concept. This is done using a script that iteratively calls ACADO for solving the CCTOC problem (1) with the non-convex inequality functions replaced by the homotopy map. The practical difference compared to Algorithm 2 is that ACADO needs to redo all the initialization when an SQP-sequence has terminated. Hence, with the scripted version we introduce a lot of overhead, and we cannot exploit information such as Hessian approximations and gradient calculations in subsequent homotopy iterations. However, the results indicate what is possible to achieve using the proposed method. For all cases presented, we have used a constant homotopy step size $\Delta_\gamma$. The SQP iterations are terminated after achieving a KKT accuracy of $10^{-2}$ for the intermediate homotopy iterations, and a final KKT accuracy of $10^{-5}$.

The results from our scripted version are compared with results from standard ACADO, as well as results obtained from the Open Motion Planning Library (OMPL) [27]. OMPL is a widely used open-source library which consists of several state-of-the-art sampled-based motion planning algorithms. In this paper, the SST planner [4] is used, since it is the only control-based planner in OMPL that is able to optimize the trajectory according to a specified objective function. All simulations were performed on an Intel Core i7-5600U processor.

### A. Motion Models

The numerical results presented in this section are based on two different types of motion models. The first one is derived from the simple car model in [1]

$$\dot{\mathbf{x}}(t) = \begin{bmatrix} \dot{x}(t) \\ \dot{y}(t) \\ \dot{\theta}(t) \\ \dot{v}(t) \\ \dot{\psi}(t) \end{bmatrix} = f(\mathbf{x}(t), \mathbf{u}(t)) = \begin{bmatrix} v(t)\cos\theta(t) \\ v(t)\sin\theta(t) \\ v(t)\tan\psi(t) \\ u_v(t) \\ u_\psi(t) \end{bmatrix}. \quad (13)$$

Here, $\mathbf{x} = [x, y, \theta, v, \psi]$ is the state vector of the car which represents the position, heading, speed and steering angle, respectively. The control signal to the system is $\mathbf{u} = [u_v, u_\psi]$ which represents the change rate in speed and steering angle, respectively. The steering angle is limited to $-\frac{\pi}{4} \leq \psi(t) \leq \frac{\pi}{4}$ and the control signals are limited to $-2 \leq u_v(t) \leq 2$ and $-\frac{\pi}{3} \leq u_\psi(t) \leq \frac{\pi}{3}$. The speed of the car can be both positive (forward motion) and negative (backward motion), limited to $-1 \leq v(t) \leq 1$.

The second model is an extension to allow movements in 3D. For this, we have used a model representing a simple fixed-wing aircraft

$$\dot{\mathbf{x}}(t) = \begin{bmatrix} \dot{x}(t) \\ \dot{y}(t) \\ \dot{\theta}(t) \\ \dot{z}(t) \\ \dot{v}(t) \\ \dot{\psi}(t) \\ \dot{\phi}(t) \end{bmatrix} = \begin{bmatrix} v(t)\cos\theta(t)\sin\phi(t) \\ v(t)\sin\theta(t)\sin\phi(t) \\ v(t)\tan\psi(t) \\ v(t)\cos\phi(t) \\ u_v(t) \\ u_\psi(t) \\ u_\phi(t) \end{bmatrix}. \quad (14)$$

This model also contains states that represent the altitude, $z(t)$, and elevation angle, $\phi(t)$, which is limited to $0 \leq \phi(t) \leq \pi$. The magnitude and direction of the velocity vector is controlled by changing $\mathbf{u} = [u_v, u_\psi, u_\phi]$. The added control signal for controlling $\phi(t)$ is limited to $-\frac{\pi}{3} \leq u_\phi \leq \frac{\pi}{3}$. For this case, we require a speed that is strictly positive, limited to $0.2 \leq v(t) \leq 1$. Both these models have been implemented in OMPL by using the available *CompoundStateSpace* class.

The objective function in the CCTOC problem will be, for both models, to minimize the total path length. This can be achieved using the objective function $f_0(\mathbf{x}(t), \mathbf{u}(t)) = |v(t)|$ in (1). In ACADO, the absolute value is formulated using a slack variable and an epigraph formulation [26]. In OMPL, the objective function is approximated with $\sum_i |v|_i dt_i$, where $i$ is the $i$th state in the solution, and $dt_i$ the amount of time $i$th control signal is applied (set to vary randomly between 0.1, 0.2, 0.3, 0.4 and 0.5 seconds).

### B. Simulation Results

In this section, solutions from the following four trajectory planning problems will be presented

- *P1:* $x_{t_0} = [1, 1, 0, 0, 0]$, $x_{t_f} = [9, 9, 0, 0, 0]$, using the motion model (13) and environment specified by (a) in Figure 4.
  SHQP: $\Delta_\gamma = 0.02$, prediction horizon $N = 30$.
  SST: Pruning radius $\delta_s$ : 0.3 , Selection Radius $\delta_v$: 1.2
- *P2:* $x_{t_0} = \left[1, 1, \frac{\pi}{2}, 1, 0\right]$, $x_{t_f} = \left[15.5, 3.5, \frac{\pi}{2}, 0, 0\right]$, using the motion model (13) and environment specified by (b) in Figure 4.
  SHQP: $\Delta_\gamma = 0.01$, prediction horizon $N = 40$.
  SST: Pruning radius $\delta_s$ : 0.3 , Selection Radius $\delta_v$: 0.6

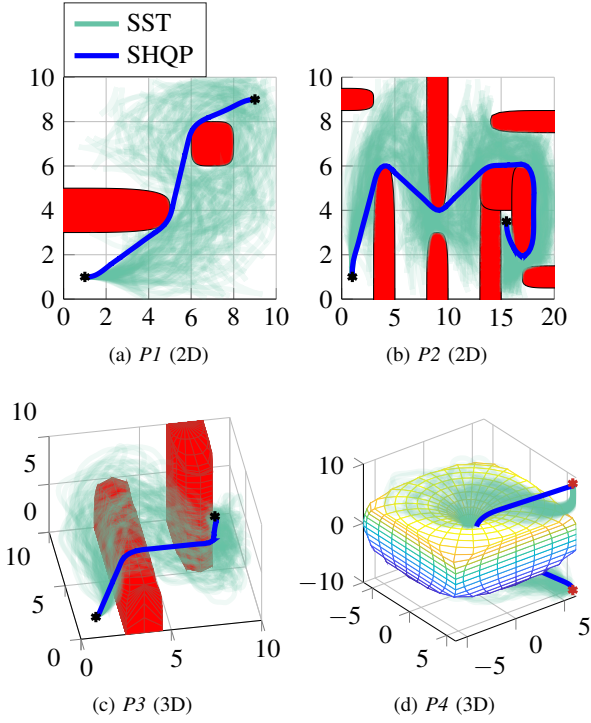

Fig. 4: The resulting trajectories for the four trajectory planning problems after 100 simulations. The variance in the solutions from SST is illustrated with opacity; higher color intensity indicates more frequently visited areas. The only way to get from initial to final state in *P4* is via the hole in the middle, which has a minimum radius of 0.5 meter.

TABLE I: The number of executions, out of 100, that returned a solution within 1000 seconds with a final state in the goal region. Standard ACADO was not able to produce a feasible solution to any of the problems.

| Planner | P1 | P2 | P3 | P4 |
|---|---|---|---|---|
| ACADO | - | - | - | - |
| SHQP | 100 | 100 | 100 | 100 |
| SST | 100 | 77 | 89 | 13 |

- *P3*: $x_{t_0} = [1, 1, \frac{\pi}{4}, 1, 1, 0, 0]$, $x_{t_f} = [9, 9, \frac{\pi}{2}, 1, 1, 0, 0]$, using the motion model (14) and environment specified by (c) in Figure 4.
  SHQP: $\Delta_\gamma = 0.02$, prediction horizon $N = 40$.
  SST: Pruning radius $\delta_s$ : 0.3 , Selection Radius $\delta_v$: 0.9
- *P4*: $x_{t_0} = [6, 6, 0, 9, 1, 0, \pi]$, $x_{t_f} = [6, 6, 0, -9, 1, 0, \pi]$, using the motion model (14) and environment specified by (d) in Figure 4.
  SHQP: $\Delta_\gamma = 0.01$, prediction horizon $N = 40$.
  SST: Pruning radius $\delta_s$ : 0.25 , Selection Radius $\delta_v$: 0.5

In all examples, the area outside what is shown in Figure 4 is prohibited. For more information regarding the parameters $\delta_s$ and $\delta_v$, the reader is referred to [4]. We terminate the SST planner when a solution is found that satisfies $||\mathbf{x}(t_f) - x_{t_f}|| < 1$, or after 1000 seconds, and returns the solution with $\mathbf{x}(t_f)$ closest to the goal state.

In Figure 4 (a) we have a trajectory planning problem with two homotopy classes, where our algorithm finds a solution in the homotopy class to the left of the Disconnected obstacle. Figure 4 (b) is a challenging maze-type problem, which requires the largest amount of computation time from our method (see Table II) due to the complex environment. The resulting computation times for *P3* and *P4* indicate that our method scales well in higher dimensions, which is one of the main advantages of the proposed method. Another advantage, compared to sampled-based planners, is that the goal state can be reached with a much higher accuracy (see Table IV) without dependencies to other methods such as steering functions.

TABLE II: Mean ($\mu_t$) and standard deviation ($\sigma_t$) of the computation time (seconds) from 100 executions of the different planners, with a maximum allowed execution time of 1000 seconds.

| Planner | P1 | P2 | P3 | P4 |
|---|---|---|---|---|
| SHQP | $\mu_t$: 3.13 $\sigma_t$: 0.02 | $\mu_t$: 28.5 $\sigma_t$: 0.23 | $\mu_t$: 9.98 $\sigma_t$: 0.06 | $\mu_t$: 7.64 $\sigma_t$: 0.04 |
| SST | $\mu_t$: 42.4 $\sigma_t$: 89.2 | $\mu_t$: 560 $\sigma_t$: 313 | $\mu_t$: 320 $\sigma_t$: 336 | $\mu_t$: 927 $\sigma_t$: 222 |

TABLE III: Mean ($\mu_o$) and standard deviation ($\sigma_o$) of the objective value (path length in [m]) from the executions that returned a solution with goal state within the goal region (see Table I).

| Planner | P1 | P2 | P3 | P4 |
|---|---|---|---|---|
| SHQP | $\mu_o$: 12.3 $\sigma_o$: 0 | $\mu_o$: 27.5 $\sigma_o$: 0 | $\mu_o$: 15.3 $\sigma_o$: 0 | $\mu_o$: 25.1 $\sigma_o$: 0 |
| SST | $\mu_o$: 20.1 $\sigma_o$: 4.40 | $\mu_o$: 57.4 $\sigma_o$: 7.07 | $\mu_o$: 24.9 $\sigma_o$: 4.47 | $\mu_o$: 38.0 $\sigma_o$: 7.14 |

One problem which is challenging for the sampled-based methods is the torus example in Figure 4 (d), since it is required to draw samples which brings the platform through the hole to get pass the obstacle [1]. This can be seen in Table I; only 13 % of the executions were able to find a solution that satisfies $||\mathbf{x}(t_f) - x_{t_f}|| < 1$ within 1000 seconds. It can also be seen that the SST method gets greedy and tries to sample straight towards the goal, while our method steers directly to the center of the torus and obtains a 100 % success rate.

From the resulting trajectories in Figure 4, it can be seen that our emulated version of Algorithm 2 manages to find a solution for all four problems and the solutions are practically relevant, with a lower objective value than the solutions produced by SST (see Table III). Also, the variance in the solutions for SST indicates that the worst case performance for this method is poor. Standard ACADO, without our proposed homotopy approach, was not able to solve any of these problems due to the non-convex inequality constraints representing the obstacles.

From our experiments, we note that it is highly relevant to use QP solvers that efficiently can solve problems with long prediction horizons $N$, such as those presented in [29], [30], [31], which have a computational complexity that grows as $\mathcal{O}(N)$. For challenging problems, the recent result in [32],

TABLE IV: Mean ($\mu_o$) and standard deviation ($\sigma_d$) of $||\mathbf{x}(t_f) - x_{t_f}||$, without considering obstacles, from 100 executions of the different planners.

| Planner | P1 | P2 | P3 | P4 |
|---|---|---|---|---|
| SHQP | $\mu_d$: 3.35e$^{-14}$ | $\mu_d$: 7.78e$^{-14}$ | $\mu_d$: 3.56e$^{-14}$ | $\mu_d$: 3.35e$^{-14}$ |
|  | $\sigma_d$: 0 | $\sigma_d$: 0 | $\sigma_d$: 0 | $\sigma_d$: 0 |
| SST | $\mu_d$: 0.79 | $\mu_d$: 1.01 | $\mu_d$: 0.86 | $\mu_d$: 1.30 |
|  | $\sigma_d$: 0.11 | $\sigma_d$: 0.57 | $\sigma_d$: 0.15 | $\sigma_d$: 0.37 |

which gives a computational complexity as low as $\mathcal{O}(\log N)$ would be useful.

## VIII. CONCLUSIONS AND FUTURE WORK

In this paper, it is shown that by combining the homotopy method and the SQP method, a method that has promising properties to solve motion planning problems involving non-convex environments is obtained. One use of the method is to compute trajectories offline. However, because of the good warm-start properties of SQP methods in general, it is also expected to be well-suited for reactive planning online when information about the environment is constantly updated. A proof-of-concept implementation of the algorithm is applied to four motion planning problems, including a challenging maze problem in 2D and two 3D problems. The result is that it handles problems which could previously not be solved by a commonly used numerical optimal control solver. Furthermore, the method outperforms an open-source sampled-based motion planner known as state-of-the-art.

The results presented in this paper can be extended in several ways. One extension is to integrate and implement Algorithm 2 in a state-of-the-art SQP solver suited for CCTOC problems with long prediction horizons. Another important future work is to develop an obstacle object and classification framework, and to automatically select suitable homotopy maps for the obstacles (*i.e.*, how the obstacles are "introduced"). Finally, it is interesting to investigate how to adaptively update the step size of the homotopy parameter more generally and flexibly by utilizing information from the SQP iterations, similar to what is done in [13] and [15].